\documentclass[a4paper,leqno,10pt,final]{amsproc}

\usepackage{amssymb,latexsym,amsxtra,amscd,amsfonts,amsthm,amsmath,verbatim}

\numberwithin{equation}{section}

\theoremstyle{plain}

\theoremstyle{definition}

\theoremstyle{remark}

\newcommand{\N}{\text{\bf N}}           

\newcommand{\Z}{\text{\bf Z}}           

\newcommand{\be}{\begin{equation}}

\newcommand{\beu}{\begin{equation*}}

\begin{document}

\title{ The K-group of Substitutional Systems}
\author{ A. El Kacimi}
\address{Lamath,  Le Mont Houy
Universit\'e de Valenciennes
59313 Valenciennes Cedex 9 (France)}
\email{aziz.elkacimi@univ-valenciennes.fr }
\author{R. Parthasarathy}
\address{ School of Mathematics, Tata Institute of Fundamental Research, Homi Bhabha Road, Colaba, Mumbai 400 005, India }
\email{sarathy@math.tifr.res.in }

\begin{abstract}
In another article we associated a dynamical system to a
non-properly ordered Bratteli diagram. In this article we describe
how to compute the $K-$group $K_0$ of the dynamical system in
terms of the Bratteli diagram. In the case of properly ordered
Bratteli diagrams this description coincides with what is already
known, namely the so-called dimension group of the Bratteli
diagram. The new ordered group defined here is more relevant for
non-properly ordered Bratteli diagrams. We use our main result to
describe $K_0$ of a substitutional system.
\end{abstract}
\maketitle

\vskip0.5cm

\noindent {\bf 0. Introduction}
\medskip

\noindent  An important tool in the study of Cantor minimal
dynamical systems $(X,T)$ is its $K$-theory; in particular the
$K_0-$group $K^0(X,T)$, which is an ordered group, is an important
invariant. After the celebrated Vershik-Herman-Putnam-Skau
approach of codifying minimal Cantor dynamical systems by using
the so-called ordered Bratteli diagrams, it became relevant to
understand the group $K_0$ directly through diagrams. This is
achieved in [HPS, Thm.5.4 and Cor.6.3] when properly ordered
Bratteli diagrams are employed. Recently, we showed how to
associate dynamical systems to non-properly ordered Bratteli
diagrams. We generalise the above result of [HPS] by a careful
modification (see 3.1) of the notion of dimension group of an ordered
Bratteli diagram. In doing this we have employed the ``tripling"
construction that was first introduced in [EP]. The result which describes the
group $K_0$ in the case of a substitutional system arising from a primitive aperiodic non-proper substitution
is described in theorem 3.12. It may be remarked that a method of computing
$K_0$ even in the case of non-proper substitutions is indicated in [DHS, sections 5,6,7]; it relies on 
showing  that the substitutional dynamical system is isomorphic to another one arising from a proper substitution.
The proof of [DHS, proposition 20] and [DHS, proposition 23]
relies heavily on `return words' and `derivative sequences' ({\it loc.cit}). But this method
seems to us to be quite indirect and not entirely transparent; `return words' are essentially of an 
existential nature and  hence do not afford a feasible method by which to compute efffectively the dimension groups
or even the Bratteli diagrams of the preferred proper substitutional system.
In contrast we feel
that our description in theorem 3.12 for non-proper substitutions is direct and 
closer in its approach and simplicity to the above cited Herman-Putnam-Skau description for
proper substitutions. It eliminates the handicap of having to first work out details in the
properly ordered case and then do the job of reducing to one such. Our methods
have the advantage of standing up equally well for the task of computing the dimension group of the dynamical system
associated to any simple non-stationary non-properly ordered Bratteli diagram (theorem 3.9).

\medskip

\noindent {\bf 1. Preliminaries}

\medskip

Since our description of the modification of {\it dimension group}
in the case of non-properly ordered Bratteli diagrams depends
heavily on the key constructions that were first introduced in
[EP] we summarize the same for the benefit of the reader following
closely the text of the first chapter of [EP]. Some of the basic
definitions and concepts in the study of Cantor dynamical systems
are also recalled in this section.

\medskip

A {\it topological dynamical system} is a pair $(X, \varphi)$ where $X$
is a compact metric space and $\varphi$ is a homeomorphism in $X$.  We say
that $\varphi$ is {\it minimal} if for any $x \in X$, the $\varphi$-orbit
of $x:= \{ \varphi^n (x)\mid  n \in \Z \}$ is dense in $X$.  We say that
$(X, \varphi)$ is a {\it Cantor dynamical system} if
$X$ is a Cantor set, i.e. $X$ is totally disconnected without isolated points.
$(X, \varphi)$ is a {\it Cantor minimal dynamical system} if, in addition,  $\varphi$ is minimal.
Some of the basic concepts of the theory are recalled below, mostly
from the more detailed sources [DHS] and [HPS].

\paragraph*{\bf 1.1. Bratteli diagram. }  A Bratteli {\it diagram} is an infinite directed
graph $(V,E)$, where $V$ is the vertex set and $E$ is the edge set.  Both
$V$ and $E$ are partitioned into non-empty disjoint finite sets
$$V= V_0 \cup V_1 \cup V_2 \cdots \ {\rm and} \ E=E_1 \cup  E_2 \cup  \cdots $$
There are two maps $r,s:E \to V$ the {\it range} and {\it source} maps.  The following
properties hold:

\begin{enumerate}
\item[(i)] $V_0 = \{v_0\}$ consists of a single point, referred to as the
`top vertex' of the Bratteli diagram
\item[(ii)] $r (E_n) \subseteq V_n, s (E_n) \subseteq V_{n-1}, n=1,2, \cdots
$.  Also $s^{-1} (v) \neq \phi$ $\forall v \in V$ and $r^{-1} (v) \neq
\phi$  for all $v \in V_1, V_2,\cdots$.
\end{enumerate}

Maps between Bratteli diagrams are assumed to preserve gradings and
intertwine the range and source maps.
If $v \in V_n$ and $w \in V_m$, where $m >n$, then a path from $v$ to
$w$ is a sequence of edges $(e_{n+1}, \cdots, e_m)$ such that $s (e_{n+1})
=v, r (e_m)= w$ and $s (e_{j+1})=r (e_j)$. Infinite paths from $v_0 \in V_0$
are defined similarly.  The Bratteli diagram is called {\it simple} if for any $n=0,1,2, \cdots
$ there exists $m>n$ such that every vertex of $V_n$ can be joined to every vertex
of $V_m$ by a path.

\paragraph*{\bf 1.2. Ordered Bratteli diagram. }  An {\it ordered} Bratteli diagram $(V, E, \geq)$ is a Bratteli
diagram $(V,E)$ together with a linear order on $r^{-1} (v), \forall v \in
V - \{v_0\} = V_1 \cup  V_2 \cup V_3\cdots $.
We say that an edge $e \in E_n$ is a {\it maximal} edge (resp. {\it minimal} edge)  if $e$ is
maximal (resp. minimal) with respect to the linear order in $r^{-1} (r (e))$.

Given $v \in V_n$, it is easy to see that there exists a unique path
$(e_1, e_2, \cdots ,e_n)$ from $v_0 $ to $v$ such that each $e_i$ is
maximal (resp. minimal).

Note that if $m >n$, then for any $w \in V_m$, the set of paths starting from
$V_n$ and ending at $w$ obtains an induced (lexicographic) linear order:
$$(e_{n+1}, e_{n+2}, \cdots, e_m) > (f_{n+1}, f_{n+2}, \cdots, f_m)$$
if for some $i$ with $n+1 \leq i \leq m, e_j=f_j$ for $ 1 <j \leq m$ and $e_i >f_i$.

\paragraph*{\bf 1.3. Proper order. }  A {\it properly ordered} Bratteli diagram is a simple ordered Bratteli diagram
$(V, E \geq )$ which possesses a unique infinite path $x_{\rm max} = (e_1,
e_2, \cdots)$ such that each $e_i$ is a maximal edge and a unique
infinite path $x_{\rm min} = (f_1, f_2, \cdots)$ such that each $f_i$ is a minimal edge.

Given a properly ordered Bratteli diagram $B=(V,E, \geq)$ we denote by $X_B$ its
infinite path space.  So
$$X_B = \{(e_1, e_2, \cdots) \mid e_i \in E_i, r (e_i) = s (e_{i+1}), i= 1,2,
\cdots \}$$
For an initial segment $(e_1, e_2, \cdots,e_n)$ we define the cylinder sets
$$U(e_1, e_2, \cdots e_n) =\{ (f_1, f_2, \cdots ) \in X_B
\mid f_i = e_i, 1 \leq i \leq n \} .$$
By taking cylinder sets to be a basis
for open sets $X_B$ becomes a topological space.  We exclude trivial cases
(where $X_B$ is finite, or has isolated points).  Thus, $X_B$ is a Cantor set.
$X_B$ is a
metric space, where for two paths $x,y$ whose initial segments to level $m$
agree but not to level $m+1, d (x,y)= 1/m+1$.

\paragraph*{\bf 1.4. Vershik map for a properly ordered Bratteli
diagram. } If $x= (e_1, e_2, \cdots e_n, \cdots) \in X_B$ and if
at least one $e_i$ is not maximal define
$$V_B (x) =y
=(f_1, f_2, \cdots, f_j, e_{j+1}, e_{j+2}, \cdots ) \in X_B$$
where $e_1, e_2, \cdots , e_{j-1}$ are maximal, $e_j$ is not maximal and has $f_j$ as successor
in the linearly ordered set $r^{-1} (r (e_j))$ and $(f_1,
f_2, \cdots, f_{j-1})$ is the minimal path from $v_0$ to $s (f_j)$.
Extend the above $V_B$ to all of $X_B$ by setting $V_B (x_{\rm max}) = x_{\rm min}$.  Then $(X_B, V_B)$ is a
Cantor minimal dynamical system.

 Next, we describe the construction of a
dynamical system associated to a non-properly ordered Bratteli diagram.
The Bratteli diagram need not be simple.
To motivate this construction, it is perhaps worthwhile to begin by
indicating how it works in the case of an ordered Bratteli diagram
associated to a nested sequence of Kakutani-Rohlin partitions of a Cantor dynamical
system $(X,T)$.

\paragraph*{\bf 1.5. K-R partition. }  A {\it Kakutani-Rohlin  partition} of the
Cantor minimal system $(X,T)$ is a clopen partition ${\mathcal P}$ of the kind
$${\mathcal P}= \{T^j Z_k \mid k \in A \ {\rm and} \ 0 \leq j <h_k\}$$
where $A$ is a finite set and $h_k$ is a positive integer.  The $k^{th}$
{\it tower} ${\mathcal S}_k$ of ${\mathcal P}$ is $\{T^j Z_k \mid  0 \leq j <h_k\}$ ; its {\it floors} are
$T^j Z_k, (0 \leq j <h_k)$.
The {\it base} of
${\mathcal P}$ is the set $Z=\bigcup_{k \in A} Z_k$.

Let  $\{{\mathcal P}_n\}, (n \in \N)$ be a sequence of Kakutani-Rohlin partitions
$${\mathcal P}_n =\{ T^j Z_{n,k} \mid k \in A_n,\hbox{and }
0 \leq j < h_{n,k} \} ,$$
with ${\mathcal P}_0 =\{X\}$
and with base $Z_n = \bigcup_{k \in A_n} Z_{n,k}$.  We say that this sequence is
{\it nested} if, for each $n$,
\begin{enumerate}
\item[(i)] $Z_{n+1} \subseteq Z_n$
\item[(ii)]$ {\mathcal P}_{n+1}$ refines the partition ${\mathcal P}_n$.
\end{enumerate}
For the Bratteli-Vershik system $(X_B, V_B)$ of sections 1.3-1.4, one obtains a Kakutani-Rohlin
partition ${\mathcal P}_n$ for each $n$ by taking the sets in the partition to be the
cylinder sets $U(e_1, e_2, \cdots e_n)$ of section 1.3 and taking as the base of the
partition the union $\bigcup U(e_1, e_2, \cdots e_n)$ over minimal paths (i.e.,
each $e_i$ is a minimal edge). This is a nested sequence.

\noindent {\bf 1.6.} To any nested sequence
$\{{\mathcal P}_n\}, (n \in \N)$ of Kakutani-Rohlin partitions we associate an ordered Bratteli diagram $B=(V,E, \geq)$
as follows (see [DHS, section 2.3]): the $\mid A_n \mid$ towers in ${\mathcal P}_n$ are in $1-1$ correspondence
with $V_n$, the set of vertices at level $n$.  Let
$v_{n,k} \in V_n$ correspond to the tower ${\mathcal S}_{n,k} =\{T^j Z_{n,k}
\mid0 \leq j < h_{n,k}\}$ in ${\mathcal P}_n$.  We refer to $T^j Z_{n,k}, 0 \leq j <
h_{n,k}$ as {\it floors} of the tower ${\mathcal S}_{n,k}$ and to $h_{n,k}$ as the {\it
height} of the tower.  {\it We will exclude nested sequences of K-R partitions where the infimum (over $k$ for fixed $n$) of
the height $h_{n,k}$ does not go to infinity with $n$.}
Let us view the tower ${\mathcal S}_{n,k}$ against the partition ${\mathcal P}_{n-1} =\{T^j Z_{n-1,k} \mid k \in A_{n-1}$, and $0 \leq j < h_{n-1,k}\}$.
As the floors of ${\mathcal S}_{n,k}$ rise from level $j=0$ to level $j=h_{n,k}-1$,  ${\mathcal S}_{n,k}$ will
start traversing a tower ${\mathcal S}_{n-1, i_1}$ from the bottom to the top floor,
then another tower ${\mathcal S}_{n-1, i_2}$ from the bottom to the top floor, then
another tower ${\mathcal S}_{n-1,i_3}$ likewise and so on till a final segment
of ${\mathcal S}_{n,k}$ traverses a tower ${\mathcal S}_{n-1, i_m}$ from the bottom to the top.  Note
that in this final step the top floor $T^j Z_{n,k}$ for $j=h_{n,k}-1$
of ${\mathcal S}_{n,k}$ reaches the top floor $T^q Z_{n-1,i_m}$ for $q= h_{n-1, i_m} -1
$ of ${\mathcal S}_{n-1, i_m}$ as a consequence of the assumption $Z_n\subset
Z_{n-1}$ and the fact that $T^{-1}$ (union of bottom floors) = union of top
floors.  Bearing in mind this order in which ${\mathcal S}_{n,k}$ traverses
${\mathcal S}_{n-1, i_1}, {\mathcal S}_{n-1, i_2}, \cdots, {\mathcal S}_{n-1, i_m}$ we associate $m$ edges,
ordered as $e_{1,k} < e_{2,k} < \cdots < e_{m,k}$ and we set the range and
source maps for edges by $r (e_{j,k})= v_{n,k}$ and $s (e_{j,k}) =v_{n-1, i_j}$.
Note that $m$ depends on the index $k \in A_n $ (and that by convention the
indexing sets $A_n$  are disjoint).  $E_n$ is the disjoint union over
$k \in A_n$ of the edges having range  in $V_n$.

\noindent {\bf 1.7.} For $x \in X$, we define $x_n \in {\mathcal P}_n^{\Z}, n \in \N$ as follows:  $x_n=
(x_{n,i})_{i \in \Z}$, where $x_{n,i} \in {\mathcal P}_n$ is the unique floor in
${\mathcal P}_n$ to which $T^i (x)$ belongs.  If $m >n$, let $j_{m,n} : {\mathcal P}_m \to {\mathcal P}_n$
be the unique map defined by $j_{m,n} (F) =F' $ if $F \subseteq F'$.  (By
abuse of notation, we use the same symbol  $F$ to denote a point of the finite
set ${\mathcal P}_m$ and also to denote the subset of $X$, in the partition ${\mathcal P}_n$,
which $F$ represents).  An important property of the map
$$X \to \prod_n
({\mathcal P}_n^{\Z}), x \mapsto (x_1, x_2, \cdots), x_n = (x_{n,i})_{i \in \Z},$$
defined above is the following:

\noindent {\bf 1.8.} If $F$ and $TF$ are two successive
floors of a ${\mathcal P}_n$-tower and if $x_{n,i}=F$ then $x_{n,i+1}=TF$. If $x_{n,i}$
 is the top floor of a ${\mathcal P}_n$-tower, then $x_{n,i+1}$ is the bottom floor of a ${\mathcal P}_n$-tower. More importantly, given
integers $K$ and $n$, there exist $m>n$ and
a single tower ${\mathcal S}_{m,k}$ of level $m$ such that the finite sequence $(x_{n,i}
)_{-K \leq i\leq  K}$ is  an interval segment contained in
$$\{j_{m,n}(T^\ell
(Z_{m,k})) \mid 0 \leq  \ell < h_{m,k} \}.$$This is a consequence of the assumption that 
the infimum of the heights of level-$n$ towers goes to infinity. 
It is true that $x_{n,i} =j_{m,n}
(x_{m,i})$, but the sequence $(x_{m,i})_{-K \leq i \leq K}$ need not be an
interval segment of $\{T^{\ell} (Z_{m,k}) \mid 0 \leq \ell <h_{m,k}\}$.

The foregoing observations in the case of an ordered Bratteli diagram associated to a nested sequence
of Kakutani-Rohlin partitions gives us the hint to define a dynamical system  $(X_B,T_B)$ of a non properly ordered
Bratteli diagram $B = (V, E, \geq )$ as follows:

\paragraph*{\bf 1.9. Definition.} For each $n$ define ${\varpi}_n =$ the set of paths from $V_0 $ to $V_n$.  There is an obvious
truncation map $j_{m,n}:{\varpi}_m \to {\varpi}_n$ which truncates paths from $V_0$
to $V_m$ to the initial segment ending in $V_n$. For
each $v \in V_n$, the set $\varpi (v)$ of paths from $\{ \ast \} \in V_0$ ending at $v$ will be called a
`${\varpi}_n$-{\it tower parametrised by} $v$'.  Each tower is a linearly ordered
set (whose elements may be referred to as floors of the tower) since
paths from $v_0$ to $v$ acquire a linear order ({\it cf.} 1.2).{\it We will exclude unusual examples
of ordered Bratteli diagram where the infimum of the height of level-$n$ towers does not go to infinity,
with $n$} (for example like [HPS, Example 3.2]).
Now, we define

\noindent{\bf 1.10.  Definition.} $X_B = \{x = (x_1, x_2, \cdots, x_n, \cdots ) \}$ where
\begin{enumerate}
\item[(i)] $x_n = (x_{n,i})_{i \in \Z} \in {\varpi}_n^{\Z}$,
\item[(ii)]  $j_{m,n} (x_{m,i}) =x_{n,i}$ for $m>n$ and $i \in \Z$    and
\item[(iii)] given $n$ and
$K$ there exists $m$ such that  $m >n$ and a vertex $v \in V_m$,
such that the interval segment $x_n [-K, K]:= (x_{n,-K},
x_{n,-K+1}, \cdots, x_{n,K})$ is obtained by applying $j_{m,n}$ to
an interval segment of the linearly ordered set of paths from
$v_0$ to $v$.
\end{enumerate}
The condition (iii) is the crucial part of the definition. Without it what one
gets is an inverse system.

The condition (iii) implies that a property similar to (1.8) holds.
Since each ${\varpi}_n$ is a finite set ${\varpi}_n^{\Z}$ has a product topology which
makes it a compact set - in fact a Cantor set.  Likewise, $\prod_n ({\varpi}_n^{\Z})$
is again a Cantor set.  Thus, $X_B \subseteq \prod_n ({\varpi}_n^{\Z})$ has
an induced topology. The lemma below and the following proposition are
analogous to corresponding facts for the Vershik model associated to properly
ordered Bratteli diagrams.

\smallskip
The following results (1.11) and (1.12) are proved in [EP].
\paragraph*{\bf 1.11. Lemma.}  {\it The topological space $X_B$ is compact.}
\smallskip

Denote by $T_B$ the restriction of the shift operator to $X_B$.
So, if $ x = (x_1, x_2, \cdots, x_n, \cdots ), $ where $x_n=(x_{n,i})_{i \in \Z} \in {\varpi}_n^{\Z}$,
then $T_B(x) = (x_1',x_2',\cdots,x_n',\cdots) $,
where $x_n'=  (x_{n,i}')_{i \in \Z} \in {\varpi}_n^{\Z}$ and $x_{n,i}'=x_{n,i+1}$.

$(X_B, T_B)
$ will be called the dynamical system associated to $B=(V, E, \geq)$.

\smallskip
\paragraph*{\bf 1.12. Proposition.} {\it If $B=(V, E, \geq)$ is a simple ordered
Bratteli diagram, then $(X_B, T_B)$ is a Cantor minimal dynamical system.}
\smallskip

\paragraph*{\bf 1.13.} In (1.7), given a nested sequence of Kakutani-Rohlin partitions of
$(X,T)$, we defined a map from $(X,T)$ to the dynamical
system $(X_B, T_B)$ of the associated ordered Bratteli diagram.  It follows
that if $(X,T)$ is minimal, and if the Bratteli diagram of the nested
sequence of K-R partitions is a simple Bratteli diagram, then
$(X,T) \to (X_B, T_B)$ is onto. If the topology of $(X,T)$ is spanned
by the collection of the clopen sets belonging to the K-R partitions
then clearly the map $(X,T) \to (X_B, T_B)$ is injective. In particular,
if the Bratteli diagram is properly ordered then the Bratteli-Vershik
system is naturally isomorphic to the system given by our construction
in 1.10.

\smallskip

\paragraph*{\bf 1.14.} Note that the same term {\it `towers'} has been used to denote two
different but related objects [in (1.5) and (1.9)]. For $v \in V_n
$, let $y$ be a path from $\{ *\}$ to $v$ in $(V, E, \geq)$.  So,
$y$ is a  `floor' (consisting of the single element $y$) belonging
to the ${\varpi}_n$- tower $\varpi (v)$ (a finite set)
parametrized by $v\in V_n$ - all in the sense of  $(1.9)$. Here,
$\varpi (v)$= all paths from $\{ *\}$ to $v$. Put ${\mathcal F}_y
= \{x= (x_1, x_2, \cdots, x_n, \cdots ) \in X_B  \mid x_{n,0} =y
\}$. ${\mathcal F}_y$ is a clopen set of the Cantor set $X_B$. Put
${\mathcal P}_n = \{{\mathcal F}_y \mid y \in \varpi (v), v \in
V_n\}.$ Then, in the sense of $(1.5)$ ${\mathcal P}_n$ is a K-R
partition of $X_B$ whose base is the union of $\bigcup {\mathcal
F}_y,( y \text{ minimal } \in \varpi (v), v \in V_n) $. Its towers
${S}_v$ are parametrized by $v\in V_n$: ${S}_v = \{{\mathcal
F}_y\mid  y\in \varpi (v)\}$. ${\mathcal F}_y,(y\in \varpi (v))$
are the floors of the tower ${S}_v$. (We encountered this K-R
partition earlier in the case of the Bratteli-Vershik system at
the end of 1.5). The ordered Bratteli diagram obtained from
$\{{\mathcal F}_y \mid y \in \varpi (v), v \in V_n\}$ is
$(V,E,\geq)$.

\medskip
\noindent {\bf 2. The Bratteli diagram $(V^{\mathcal Q}, E^{\mathcal Q}, \geq )$ }
\medskip

\noindent{\bf 2.1.} We will now define two nested sequences of K-R partitions of
$X$.  For $v \in V_n $,
let $y$ be a path from $\{ *\}$ to $v$ in $(V, E,
\geq)$.  So, $y$ is a  `floor'  belonging to the $\varpi_n$-
tower $\varpi (v)$ parametrized by $v$.
Put ${\mathcal F}_y = \{x= (x_1, x_2, \cdots, x_n, \cdots ) \in X  \mid
x_{n,0} =y \}$
$${\mathcal P}_n = \{{\mathcal F}_y \mid y \in \varpi(v), v \in V_n \}.$$
Then $\{{\mathcal P}_n \}_n$ is a nested sequence of K-R partitions of $X$
.  But, the topology of $X$ need not be spanned by
the collection of clopen sets $\{{\mathcal F}_y\}, (y \in \varpi (v), v
\in V_n, n \in \N)$.  In contrast, the topology of $X$
is indeed spanned by the collection of clopen sets in another nested sequence
$\{{\mathcal Q}_n\}_n$ of K-R partitions, defined below.
Let $\varpi = \varpi (u), \varpi' = \varpi (v),
\varpi'' = \varpi (w)$
be three $\varpi_{n}$-towers and $y$ a floor of $\varpi'$.
For any $x \in X$ and for any $n$ if $x_{n,i}$ is a floor of a
$\varpi_n$-tower $\overline{\varpi}$, then for some
$a,b\in \Z$ such that $a \leq i \leq b$, the segment $x_n [a,b]$ is
just the sequence of floors in $\overline{\varpi}$.  We define

${\mathcal F}(\varpi, \varpi', \varpi''; y) =$ the clopen subset of
${\mathcal F}_y$ consisting of the elements $x= (x_1, x_2, \cdots, x_n, \cdots)$
with the property that for some $a_1 < a_2 \leq 0 < a_3 < a_4 \in \Z$, the segment $x_n [a_1,
a_2-1]$ is the sequence of floors of $\varpi$, the segment $x_n
[a_2, a_3 -1]$ is  the sequence of floors of $\varpi'$ and the segment
$x_n [a_3, a_4]$ is the sequence of floors of $\varpi''$.
Some of the sets  ${\mathcal F}(\varpi, \varpi', \varpi'';y)$ may be empty, but
the non-empty
sets ${\mathcal F}(\varpi, \varpi', \varpi'';y)$ form a K-R partition which we denote by
${\mathcal Q}_n$. For fixed $\varpi, \varpi', \varpi'' $
the subcollection $\{{\mathcal F}(\varpi,
\varpi', \varpi'';y)\}$ as $y$ varies through  the floors of
$\varpi'$, is a ${\mathcal Q}_n$-tower parametrized by
$[u, v,w]$.  We denote this ${\mathcal Q}_n$-tower by
${\mathcal S}_{(\varpi,
\varpi', \varpi'')}$.  The floors of the tower ${\mathcal S}_{(\varpi, \varpi', \varpi'')}$
are $\{{\mathcal F}(\varpi, \varpi', \varpi'';y)\}$ as $y$
runs through the sequence of floors of $\varpi'$.

\noindent{\bf 2.2. The tripling of $(V, E, \geq)$.} Let
$(V, E, \geq)$ be an arbitrary simple, ordered Bratteli diagram.  Define
$(V^{\mathcal Q}, E^{\mathcal Q}, \geq )$ as follows:  $V^{\mathcal Q}_0 =\{*\}$, a
single point.

$V^{\mathcal Q}_n$ consists of triples $(u,v,w) \in V_n\times V_n \times V_n$ such that for
some $y \in V_m$ where $m>n$, the level-$m$ tower ${\varpi}(y)$ passes successively
through the level-$n$ tower ${\varpi}(u)$, then ${\varpi}(v)$ and then ${\varpi}(w)$.
An edge $ \tilde{e} \in  E^{\mathcal Q}_n$ is a triple $(u,e,w)$ such that
$e$ is an edge of $(V,E) $ and $(u,r (e),w) \in V^{\mathcal Q}_n$.  Let
\begin{enumerate}
\item[]$\{e_1, e_2, \cdots, e_k\}$ be all  the edges  in  $r^{-1} (r(e))$,
\item[]$\{f_1, f_2, \cdots, f_{\ell}\}$ be all the edges in  $r^{-1} (u)$ and
\item[]$\{g_1, g_2, \cdots, g_m \}$ be all the edges in  $r^{-1} (w)$.
\end{enumerate}
The sources of $(u,e_1, w), (u, e_2, w), \cdots, (u,e_k, w)$ are defined
to be \\
$ (s(f_{\ell}), s(e_1),s(e_2)), (s(e_1),s (e_2), s(e_3)),
\cdots, (s(e_{k-1}), s(e_k), s(g_1))$ \\
respectively. The range of
$(u,e,w)$ is of course $(u,r (e), w)$. If   $r^{-1} (r(e))$ is ordered as $\{e_1, e_2, \cdots, e_k\}$,
we declare the ordering of $r^{-1} (r(u,e,w))$ to be $\{(u,e_1, w), (u, e_2, w),\\
\cdots, (u,e_k, w)\}$. The ordered Bratteli diagram
$(V^{\mathcal Q}, E^{\mathcal Q}, \geq )$ thus defined will be called the tripling of $(V, E, \geq)$.

The map $\pi:(V^{\mathcal Q}, E^{\mathcal Q},\geq)\longrightarrow
(V,E,\geq)$ given by $(u,v,w) \mapsto v, (u,e,w) \mapsto e$ enjoys
the `{\it unique path lifting}' property in the following sense.
If $m > n \geq 1$, and $(e_n, e_{n+1}, \cdots, e_m)$ is a path in
$(V,E)$ from $V_{n-1}$ to $V_m$ with $r(e_m)=v$ then for any
$(u,v,w) \in V^ {\mathcal Q}_m$, there is a unique path
$(\tilde{e}_n, \tilde{e}_{n+1}, \cdots, \tilde{e}_{m})$ in
$(V^{\mathcal Q}, E^{\mathcal Q})$ which maps onto $(e_n, e_{n+1},
\cdots, e_m)$ under $\pi$ and such that $r (\tilde{e}_m) =
(u,v,w)$. It is quite elementary to check that the map
$\pi:(V^{\mathcal Q}, E^{\mathcal
Q},\geq)\longrightarrow(V,E,\geq)$ induces an isomorphism between
the corresponding dynamical systems given by 1.10 ([EP, 2.20]).

[ Two different edges on the left with the same source may map
into the same edge on the right. Two different edges on the left
with the same range cannot map to the same edge on the right.]

 Let $ \{n_k\}^{\infty}_{k=0}$ be a subsequence of $\{0,
1,2, \cdots \}$ where we assume $n_0=0$.  A Bratteli diagram $(V',E')$
is called a {\it `telescoping'} of $(V,E)$ if $V'_k=V_{n_k}$ and
$E'_k$ consists of paths $(e_{n_{k-1}+1}, \cdots, e_{n_k})$ from
$V_{n_{k-1}}$ to $V_{n_k}$ in $(V,E)$, the range and source maps being the obvious
ones.  It is easy to see that tripling is compatible with telescoping.

\noindent {\bf 2.3. Stationary Bratteli diagrams.} A Bratteli
diagram is stationary if the diagram repeats itself after level
$1$. (One may relax by allowing a period from some level onwards;
but, a telescoping will be stationary in the above restricted
sense.) If $(V,E,\geq )$ is an ordered Bratteli diagram and the
diagram together with the order repeats itself after level $1$,
then $(V, E, \geq)$ will be called a stationary ordered Bratteli
diagram. We refer the reader to [DHS, section (3.3)] for the usual
definition of a substitutional system and how they give rise to
stationary Bratteli diagrams. Some details are recalled below. Let
$(V,E, \geq)$ be as above and suppose moreover that it is a simple
Bratteli diagram. We have

\begin{enumerate}
\item[ (1)] an enumeration $\{v_{n,1}, v_{n,2}, \cdots, v_{n,L}\}$ of $V_n, \forall
n \geq 1$,

\item[ (2)]  for $n>1$ and $1 \leq j \leq L$ an enumeration
$\{e_{n,j,1}, e_{n,j,2}, \cdots, e_{n,j,a_j}\} $ of $\ r^{-1} (e_{n,j})$
which is assumed to be listed in the linear order in $r^{-1} (v_{n,j})$,

\item[ (3)]  in the enumerations above, $L$ does not depend on $n$
and $a_j$ depends only on $j$  and not on $n$.   Moreover, the
ordering in $\ r^{-1} (e_{n,j})$ is stationary, i.e, if $n, m >1$,
if $ 1 \leq j \leq L, 1 \leq  k \leq L, 1 \leq i \leq a_j$, then
$`` s (e_{n,j,i}) = v_{n-1,k}$'' $\Longrightarrow ``s (e_{m,j,i})
= v_{m-1, k}$''.
\end{enumerate}

\noindent{\bf 2.4. Substitutional systems.} Let $A$ be an alphabet set.
Write $A^+$ for the set of words of finite length in the alphabets of $A$.
Let $\sigma:A \to A^+ $ be a primitive aperiodic non-proper substitution, written,
$\sigma (a) = \alpha \beta \gamma \cdots $. 
The stationary ordered Bratteli diagram $B=(V,E,\geq )$ associated to
$(A, \sigma)$ ({\it cf.} [DHS, section 3.3] can be
described as
$$V_n =A , \forall \text{ }n \geq 1, V_0 =\{*\}$$
$$E_n =\{(a,k,b) \mid  a,b \in A, k \in \N,\text{ } a \text{ is the }k^{th} \text{ alphabet
in the word }\sigma  (b)\}.$$
(The reader who prefers a more carefully evolved notation can consider intoducing
an extra factor `$\times \{n\}$' so that vertices and edges at different levels are
seen to be disjoint).
The source and range maps $s$ and $r$ are
defined by $s (a,k,b)=a, r (a, k,b)=b$.  In the linear order in $r^{-1}
(b)$, $ (a,k,b)$ is the $k^{th}$ edge.
\smallskip

To the stationary ordered Bratteli diagram $B$ of $(A,\sigma)$
(which may not be properly ordered unless $\sigma$ is a primitive,
aperiodic, proper substitution, -- see [DHS, section 3]) we can
associate a dynamical system $X_B$ following the construction of
1.10; this is naturally isomorphic to the substitutional dynamical
system $(X_{\sigma},T_{\sigma})$ associated to $(A,\sigma)$ defined for example
in [DHS, section 3.3.1]. (See [EP, section 2.5].)

\noindent{\bf 2.5. Tripling for a substitutional system $(A,
\sigma)$.} Let $(A, \sigma)$ be a substitutional system and
suppose $B=(V,E,\geq )$ is the stationary ordered Bratteli diagram
associated to $(A, \sigma)$. Define $ A^{\mathcal Q} = \{(a,b,c)
\in A \times A \times A \mid abc$ {\it occurs as a subword of }
$\sigma^n (d)$ {\it for some }$d \in A $ {\it and some }$n \}$ .
Define
$$\sigma^{\mathcal Q}: A^{\mathcal Q} \to (A^{\mathcal Q})^+$$
by $\sigma^{\mathcal Q} [(a,b,c)] =(a_m, b_1, b_2) \cdot
(b_1,b_2,b_3) \cdots (b_{n-2}, b_{n-1}, b_n) \cdot (b_{n-1}, b_n,
c_1)$,  where $\sigma (b) =b_1 \cdot b_2 \cdots b_n$, and $a_m$ is
the last alphabet in $\sigma (a)$, while $c_1$ is the first
alphabet in $\sigma (c)$. Then $(V^{\mathcal Q}, E^{\mathcal Q},
\geq )$ is the stationary ordered Bratteli diagram associated to
$(A^{\mathcal Q},\sigma^{\mathcal Q})$.

\medskip
\noindent {\bf 3. The groups $K^0(X,T),K_-0(V,E,\geq)$ and $K_0(V,E)$.}
\medskip

\noindent{\bf 3.1.Definition.} Let $(X,T)$ be a Cantor minimal system.  Let $C(X,\Z)$ be the space of
integer valued continuous functions on $X$.  Let
$$
K^0(X,T) = C(X,\Z)/\partial_T C(X,\Z)
$$
where $\partial_T : C(X,\Z) \to C(X,\Z)$ denotes the coboundary operator $\partial_T(f)=f-f\circ T$.
A function of the form $f-f\circ T$ is called a coboundary. Define the positive cone
$$
K^0(X,T)^+ = \{[f] \mid f \in C(X,\Z^+)\}
$$
where $ [f]$ denotes the projection modulo coboundaries. The ordered group $(K^0(X,T),\\K^0(X,T)^+)$
has a distinguished order unit, namely [1], the projection of the constant function 1.

 Let $(V,E)$ be a Bratteli diagram and $(V,E,\geq)$ the same thing equipped with a linear order
on edges which makes it an ordered Bratteli diagram. As usual the
dimension group $K_0(V,E)$ is defined to be the inductive limit of
the system of ordered groups
$$
\Z^{|V_0|} \buildrel {A_0} \over \longrightarrow \Z^{|V_1|} \buildrel {A_1} \over \longrightarrow
\Z^{|V_2|} \buildrel {A_2} \over \longrightarrow \Z^{|V_3|} \buildrel {A_3} \over \longrightarrow
\cdots
$$
where the positive homomorphism $A_n$ is given by matrix multiplication with the incidence
matrix between levels $n-1$ and $n$. The inductive limit $K_0(V,E)$ is endowed with the induced
order, the positive cone being denoted by $K_0(V,E)^+$. The image of $1\in \Z^{|V_0|}$ in
$(K_0(V,E),K_0(V,E)^+)$ is an order unit.

 On the other hand we define the group $K_-0(V,E,\geq)$ in the following way.  Whenever we have
$m \geq n$ and two paths $\varpi_1$ and $\varpi_2,(\varpi_1 \leq
\varpi_2)$ from $V_n$ to $V_m$ with the same range $u\in V_m$
define $[\varpi_1,\varpi_2) $ to be the set consisting of all
paths from $V_n$ to $V_m$ lying between $\varpi_1 \text{
(included) and }\varpi_2 \text{ (excluded) }$ranging at $u$. Put $
B \Z^{|V_n|}=\{\overline{m}= (m_k)_{k\in V_n} \in \Z^{|V_n|}\mid
\Sigma_{{\varpi \in [\varpi_1,\varpi_2)}}  m_{s(\varpi)}=0 \text{
for all }m\text{ and all }\varpi_1,\varpi_2 $ as above with the
same source and same range$\}$.
 Observe that $A_n(B \Z^{|V_n|}) \subseteq B \Z^{|V_{n+1}|}$.
 Moreover, suppose $\overline{p},\overline{q} \in \Z^{|V_n|,+},\overline{m}\in
B\Z^{|V_n|}$ and $\overline{p} = -\overline{q} +\overline{m}$.
Then, $\overline{m}=\overline{p}+\overline{q} \in \Z^{|V_n|,+}$,
which forces $\overline{m}$ to be zero because of the defining
conditions of $B\Z^{|V_n|}$. Thus, the natural order in
$\Z^{|V_n|}$ induces an order in the quotient group
$\Z^{|V_n|}/B\Z^{|V_n|}$ making it an ordered group. Define
$(K_-0(V,E,\geq),K_-0(V,E,\geq)^+)$ to be the inductive limit of
the system of ordered groups
$$
\displaystyle\frac{\Z^{|V_0|}}{B\Z^{|V_0|}} \buildrel {A_0} \over \longrightarrow
\displaystyle\frac{\Z^{|V_1|}}{B\Z^{|V_1|}} \buildrel {A_1} \over \longrightarrow
\displaystyle\frac{\Z^{|V_2|}}{B\Z^{|V_2|}} \buildrel {A_2} \over \longrightarrow
\displaystyle\frac{\Z^{|V_3|}}{B\Z^{|V_3|}} \buildrel {A_3} \over \longrightarrow
\cdots
$$
Observe that $B\Z^{|V_0|}=0$. The image of $1\in \Z^{|V_0|}$ in
$(K_-0(V,E,\geq),K_-0(V,E,\geq)^+)$ is an order unit.

\smallskip
\noindent{\bf 3.2. Theorem.} {\it For $B=(V,E,\geq)$ let
$(X_B,T_B)$ be defined as in (1.10). Write $(X,T)=(X_B,T_B)$.
Define the tripling $B^{\mathcal Q}=(V^{\mathcal Q},E^{\mathcal
Q},\geq)$ as in 2.2.  Then $K^0(X,T)$ is naturally order
isomorphic to $K_-0( V^{\mathcal Q},E^{\mathcal Q},\geq), $
preserving distinguished order units. }

\noindent{\bf Proof.} We recall the notation introduced in 2.1. Given $f \in C(X,\Z)$, choose $n$
sufficiently large such that $f,\partial_T(f)$ are both constant on the sets of the partition $\mathcal Q_n$.
The vertices of the Bratteli diagram $(V^{\mathcal Q},E^{\mathcal Q},\geq)$ correspond to towers
${\mathcal S}_{(\varpi,\varpi', \varpi'')}$ of a K-R partition which in turn are partitioned into floors
$\{{\mathcal F}(\varpi,
\varpi', \varpi'';y)\}$ as $y$ varies through  the floors of
$\varpi'$. For $f $ as above, define $\gamma_n (f) \in \Z^{|V^{\mathcal Q}_n|}$ by
$\gamma_n (f) (\varpi,\varpi', \varpi'' ) = f(x) + f(Tx) + f(T^2x) + \cdots +f(T^{h-1}x)$ where $x$ belongs to
the lowest floor of ${\mathcal S}_{(\varpi,\varpi', \varpi'')} $ and $h$ is the height of the tower
${\mathcal S}_{(\varpi,\varpi', \varpi'')}$.
Then  $ A^{\mathcal Q}_n(\gamma_n(f)) = \gamma_{n+1}(f)$ and
$ \gamma_n( \partial_T(f)) \in B\Z^{|V^{\mathcal Q}_n|}$. This gives rise to a map
$$
\gamma : K^0(X_B,T_B) \longrightarrow K_-0(V^{\mathcal Q},E^{\mathcal Q},\geq).
$$

\smallskip
\noindent{\bf 3.3. Lemma.} {\it  Let $f \in C(X,\Z)$ and suppose
that $f$ is constant on the sets of the partition $\mathcal Q_n$.
Suppose that $\gamma_n(f) \in B\Z^{|V^{\mathcal Q}_n|}$. Then,
$f=\partial_T(g)$, for some $g\in C(X,\Z)$. }

\smallskip
\noindent{\bf 3.4. Lemma.} {\it With the same assumptions as in
3.3 suppose that $x,y (\in X) $ both lie in the same floor of a
$\mathcal Q_n$-tower ${\mathcal S}_{(\zeta,\zeta', \zeta'')}$.
Furthermore, suppose that for some positive integers $k,\ell$ both
$T^kx$ and $T^{ \ell }y$ lie in the same floor of a $\mathcal
Q_n$-tower ${\mathcal S}_{(\vartheta,\vartheta', \vartheta'')}$.
Then,
$$f(x)+f(Tx)+\cdots+f(T^kx) = f(y)+f(Ty)+\cdots+f(T^{\ell}y).$$}

\noindent{\bf Proof of 3.4.} Let $m\geq n$. Let $U_y$ be a neighborhood of $y$ such that $\forall z
\in U_y$ and for $i \in [0,\ell], T^iy$ and $T^iz$ belong to the same floor of the $\mathcal Q_n$ partition.
Since the orbit of $T^kx$ by iterations of $T$ is dense $\exists j$ such that $T^{j+k}x \in U_y$. For
sufficiently large $m\geq n, \exists $ a $\mathcal Q_m$-tower $\mathcal S$ such that $x, T^kx, T^{k+j}x,
T^{k+j+\ell}x$ belong to different floors of $\mathcal S$.

Let $u$ be the vertex of $V^{\mathcal Q}_m$ represented by
$\mathcal S$. The floors of the $\mathcal Q_m$-tower $\mathcal S$
are linearly ordered reflecting the linear order in the set of
paths in $(V^{\mathcal Q}, E^{\mathcal Q},\geq)$ from the top
vertex to $ u$. Similarly, the paths from $ V^{\mathcal Q}_n$ to
$u \in  V^{\mathcal Q}_m$ are linearly ordered reflecting the
order in which $ \mathcal S$ traverses the level-$n$ towers of
$(V^{\mathcal Q}, E^{\mathcal Q},\geq)$. Denote by
$\varpi_1,\varpi_2,\varpi_3,\cdots,\varpi_L$ the paths from $
V^{\mathcal Q}_n$ to $u \in  V^{\mathcal Q}_m$ in their linear
order. Write
$s(\varpi_1),s(\varpi_2),s(\varpi_3),\cdots,s(\varpi_L) \in
V^{\mathcal Q}_n $ for their sources and $\mathcal
S_{s(\varpi_1)},\mathcal S_{s(\varpi_2)}, \mathcal
S_{s(\varpi_3)},\cdots,\\  \mathcal S_{s(\varpi_L)}$ for the
$\mathcal Q_n$-towers represented by these sources. Thus,
$\mathcal S$  traverses $\mathcal Q_n$-towers in the order
$\mathcal S_{s(\varpi_1)},\mathcal S_{s(\varpi_2)}, \mathcal
S_{s(\varpi_3)},\cdots,  \mathcal S_{s(\varpi_L)} $. Choose $1\leq
a<b<c<d\leq L$ such that $x$,(resp.$T^kx$, resp.$T^{k+j}x$,
resp.$T^{k+j+\ell}x$), is picked up by $\mathcal S$ at the
$a^{th}$(resp.$b^{th}$, resp.$c^{th}$, resp.$d^{th}$) instance of
$\mathcal S$ traversing through a $\mathcal Q_n$-tower, namely,
$\mathcal S_{s(\varpi_a)}$,(resp.$\mathcal S_{s(\varpi_b)}$,
resp.$\mathcal S_{s(\varpi_c)}$, resp.$\mathcal S_{s(\varpi_d)}$).
In particular, observe that $s(\varpi_a) = s(\varpi_c)$ and
$s(\varpi_b) = s(\varpi_d)$.

Since $\gamma_n(f) \in B\Z^{|V^{\mathcal Q}_n|}$, we have
$$
\gamma_n(f)(s(\varpi_a)) + \gamma_n(f)(s(\varpi_{a+1})) +\cdots +\gamma_n(f)(s(\varpi_{c-1})) =0
$$ and similarly,
$$
\gamma_n(f)(s(\varpi_b)) + \gamma_n(f)(s(\varpi_{b+1})) +\cdots +\gamma_n(f)(s(\varpi_{d-1})) =0.
$$ These two equations imply that
$$
f(x)+f(Tx)+\cdots+f(T^{k+j-1}x)= 0
$$ and
$$
f(T^{k+1}x)+f(T^{k+2}x)+\cdots+f(T^{k+j+\ell}x)= 0.
$$ Hence,
$$
f(x)+f(Tx)+\cdots+f(T^kx)=
$$
$$
f(x)+f(Tx)+\cdots+f(T^kx) +f(T^{k+1}x)+\cdots+   f(T^{k+j}x)+\cdots+     f(T^{k+j+\ell}x)=
$$
$$
f(T^{k+j}x)+f(T^{k+j+1}x)+\cdots+f(T^{k+j+\ell}x)=
$$
$$
f(y)+f(Ty)+\cdots+f(T^{\ell}y).
$$
This ends the proof of 3.4.
\hfill $\square$
\smallskip

\noindent{\bf Proof of 3.3.} Choose $x_0 \in X$. Now, for any $z
\in X$ choose $k \in \Z^+$ such that $T^kx_0 $ and $z$ belong to
the same $\mathcal Q_n$-floor. Define $g \in C(X,\Z)$ by $g(z) =
f(x_0)+f(Tx_0)+\cdots+f(T^kx_0)$. Then, 3.4 implies that $g$ is
well defined and $\partial_T(g) = -f\circ T$. So,
$\partial_T(-g\circ T^{-1}) = f$.

 From 3.3 one can immediately deduce that the map $\gamma :K^0(X_B,T_B)
\longrightarrow K_-0(V^{\mathcal Q},E^{\mathcal Q},\geq)$ defined
just before the statement of lemma 3.3 is an isomorphism.

This completes the proof of Theorem 3.2.
\hfill $\square $

\noindent{\bf 3.5. A subgroup of $B\Z^{|V^{\mathcal Q}_n|}$.} In
practice it is quite tedious to determine whether a given element
$\overline{p}$ of $\Z^{|V^{\mathcal Q}_n|}$ lies in
$B\Z^{|V^{\mathcal Q}_n|}$. We now begin to describe a subgroup
$\Delta \Z^{|V^{\mathcal Q}_n|} \subseteq B\Z^{|V^{\mathcal
Q}_n|}$, which is more easily identifiable than $B\Z^{|V^{\mathcal
Q}_n|}$. Eventhough, in general, this inclusion is proper we will
later see that the distinction disappears when one takes inductive
limits. As a consequence, we are able to obtain theorem 3.9, which
yields a feasible method to compute $K_0$ effectively. Clearly,
$\Z^{|V^{\mathcal Q}_n|}$ is the space of integral valued
functions on the set $V^{\mathcal Q}_n$. For a function $\varphi
:V_n\times V_n \longrightarrow \Z$ define $\delta (\varphi) \in
\Z^{|V^{\mathcal Q}_n|}$ by $\delta (\varphi)(a,b,c) = \varphi
(b,c) - \varphi (a,b)$.

\smallskip

\noindent{\bf Lemma 3.6.} $\delta (\varphi) \in B\Z^{|V^{\mathcal
Q}_n|} $.

\smallskip

\noindent{\bf Proof.} Let $\overline{p} \in \Z^{|V^{\mathcal
Q}_n|}$. Write $\overline{p} = \{p_{(u,v,w)}\}_{(u,v,w)\in
V^{\mathcal Q}_n}$. Take two paths from $V^{\mathcal Q}_n$ to
$V^{\mathcal Q}_m(m >n)$ with the same source in $V^{\mathcal
Q}_n$ and same range in $V^{\mathcal Q}_m$. The sequence of
sources of paths lying between the above two paths is of the form
$\{(u_1,v_1,w_1), (u_2,v_2,w_2),\cdots,(u_j,v_j,w_j)\}$ where

\begin{enumerate}
 \item[ (i)] $ (u_1,v_1,w_1)=(u_j,v_j,w_j)$,
 \item[ (ii)]$u_{i+1}=v_i$ and
 \item[ (iii)] $v_{i+1}=w_i$, for $i= 1,2,\cdots,j-1$.
\end{enumerate}

If $\overline{p}=\delta(\varphi)$, the sum $p_{(u_1,v_1,w_1)}
+p_{(u_2,v_2,w_2)}+\cdots+p_{(u_{j-1 },v_{j-1 },w_{j-1 })}$
equals\\
$\{\varphi(v_1,w_1)-\varphi(u_1,v_1)\}+\{\varphi(v_2,w_2)-\varphi(u_2,v_2)\}+
\cdots +\{\varphi(v_{j-1},w_{j-1})-\varphi(u_{j-1},v_{j-1})\}$\\
$$
\begin{matrix}
 =& -\varphi(u_1,v_1)+\varphi(v_{j-1},w_{j-1})&(\text{in view of (ii) and (iii) above})\\
 =& -\varphi(u_1,v_1)+\varphi(u_j,v_j) &(\text{in view of (ii) and (iii) above})\\
 =& 0 &(\text{in view of (i) above}).
\end{matrix}
$$
The condition for $\overline{p}$ to belong to $B\Z^{|V^{\mathcal
Q}_n|}$ is precisely that the sums of the type $p_{(u_1,v_1,w_1)}
+p_{(u_2,v_2,w_2)}+\cdots+p_{(u_{j-1 },v_{j-1 },w_{j-1 })}$ as
above should all vanish. As the foregoing calculation shows this
holds whenever $\overline{p}= \delta(\varphi)$ for some $\varphi :
V_n\times V_n \longrightarrow \Z$.

\hfill $\square $

\paragraph{}  We define
$\Delta \Z^{|V^{\mathcal Q}_n|}$ to be the subgroup $\delta (\Z^{|V_n
\times V_n|})$  of $\Z^{|V^{\mathcal Q}_n|}$.

\paragraph{} Recall the map $A^{\mathcal Q}_n: \Z^{|V^{\mathcal
Q}_n|}\longrightarrow \Z^{|V^{\mathcal Q}_{n+1}|}$ given by matrix
multiplication by the incidence matrix between the levels
$V^{\mathcal Q}_n$ and $V^{\mathcal Q}_{n+1}$. For a function
$\varphi :V_n \times V_n \longrightarrow \Z$ define
$\varphi':V_{n+1} \times V_{n+1} \longrightarrow \Z$ by $\varphi'
(u',v')= \varphi (u,v)$ where $u$(resp.$v$) is the source of the
last(resp.first) edge ranging at $u'$(resp.$v'$).

\smallskip

\noindent{\bf Lemma 3.7. }{\it With notation as above,
$A^{\mathcal Q}_n (\delta (\varphi)) = \delta (\varphi')$. In
particular the identity endomorphism of $\Z^{|V^{\mathcal Q}_n |}$
induces a map
$$\displaystyle\frac{\Z^{|V^{\mathcal Q}_n|}}{\Delta\Z^{|V^{\mathcal Q}_n
|}} \buildrel {\rho_n} \over \longrightarrow
\displaystyle\frac{\Z^{|V^{\mathcal Q}_n |}}{B\Z^{|V^{\mathcal
Q}_n |}}$$ and
$$\begin{matrix}
 \displaystyle\frac{\Z^{|V^{\mathcal Q}_n|}}{B\Z^{|V^{\mathcal Q}_n|}}&\buildrel {A^{\mathcal Q}_n} \over
\longrightarrow &\displaystyle\frac{\Z^{|V^{\mathcal
Q}_{n+1}|}}{B\Z^{|V^{\mathcal Q}_{n+1}|}}   \\
&&\\
\rho_n\uparrow & & \rho_{n+1}\uparrow \\
&& \\
\displaystyle\frac{\Z^{|V^{\mathcal Q}_n|}}{\Delta\Z^{|V^{\mathcal
Q}_n|}}&\buildrel {A^{\mathcal Q}_n} \over \longrightarrow
&\displaystyle\frac{\Z^{|V^{\mathcal
Q}_{n+1}|}}{\Delta\Z^{|V^{\mathcal Q}_{n+1}|}}
\end{matrix}$$
is commutative. }

\smallskip

\noindent{\bf Proof. } The proof of the assertion $A^{\mathcal
Q}_n (\delta (\varphi)) = \delta (\varphi')$ is a straightforward
calculation using definitions and notation. The rest follows
immediately.

\hfill $\square $

 We might wish to ask whether every element
$\overline{p}$ of $ B\Z^{|V^{\mathcal Q}_n|}$ is of the form $
\delta (\varphi)$ for some function $\varphi :V_n \times V_n
\longrightarrow \Z$. The proposition 3.8 below shows that after
applying a finite iteration $A^{\mathcal Q}_{n+i}\circ \cdots
A^{\mathcal Q}_{n+1}\circ A^{\mathcal Q}_n $ to $\overline{p}$ it
will indeed be so.

\smallskip

\noindent{} Let $\overline{p} \in B\Z^{|V^{\mathcal Q}_n|}$. Let
$g \in C(X,\Z)$ be chosen as in Lemma 3.3 so that
\begin{enumerate}
 \item[ (i)] $\partial_T(g)$ is constant on the sets of the
 partition $\mathcal Q_n$ and moreover, $\overline{p} = \gamma_n (
 \partial_T(g))$.
  \item[ (ii)] $g$ itself is constant on the sets of the partition
 $\mathcal Q_{n+i}$ for some positive integer $i$.
\end{enumerate}
Choose a positive integer $j$ such that any $\mathcal
Q_{n+i+j}~$-tower traverses through at least two $\mathcal
Q_{n+i}~$-towers. Then, of course, any $\mathcal P_{n+i+j}~$-tower
traverses through at least two $\mathcal P_{n+i}~$-towers.

For $u',v'\in V_{n+i+j}$ let $u_a$(resp.$u_{a-1}$,~resp.$v_1,$
resp.$v_2$) be the source of the {\it last}(resp.{\it last but
one}, resp.{\it first}, resp.{\it second}) path from $V_{n+i}$ to
$V_{n+i+j}$ ranging at $u'$(resp.$u'$,resp.$v'$,resp.$v'$). If
there exists $x\in X$ such that
\begin{enumerate}
 \item[ (a)] $x \in $ bottom floor of the $Q_{n+i}$-tower
 represented by $(u_a,v_1,v_2)$,
 \item[ (b)] $T^{-1}x \in$ top floor of the $Q_{n+i}$-tower
 represented by $(u_{a-1},u_a,v_1)$
\end{enumerate}
define $\varphi'(u',v')=g(x)$; then, $\varphi'(u',v')$ is
independent of $x$. For given $u',v' \in V_{n+i+j}$ if no such $x$
exists, define $\varphi'(u',v')$ arbitrarily.

\smallskip

\noindent{\bf Proposition 3.8. }With notation as above
$A^{\mathcal Q}_{n+i+j-1}\circ \cdots A^{\mathcal Q}_{n+1}\circ
A^{\mathcal Q}_n ~(\overline{p}) = \delta (\varphi')$.

\smallskip

\noindent{\bf Proof. }$\delta (\phi')(u',v',w')=\phi'(v',w')
-\phi'(u',v') = g(T^hy) -g(y)$, if $y$ lies in the lowest floor of
the $\mathcal Q_{n+i+j}~$-tower of height $h$ represented by
$(u',v',w')$. Also, for $(u,v,w) \in V^{\mathcal Q}_n$, $
\overline{p}(u,v,w) = \gamma_n (\partial_T(g))(u,v,w) =$ the sum
$(g\circ T - g)(z) +(g\circ T - g)(Tz)+(g\circ T - g)(T^2z)+
\cdots +(g\circ T - g)(T^{k-1}z)$, where $k$ is the height of the
$\mathcal Q_n~$-tower represented by $(u,v,w)$ and $z$ lies in its
lowest floor. Thus, $A^{\mathcal Q}_{n+i+j-1}\circ \cdots
A^{\mathcal Q}_{n+1}\circ A^{\mathcal Q}_n
~(\overline{p})(u',v',w')$ is the sum of $g\circ T - g$ taken over
all the floors of the $\mathcal Q_{n+i+j}~$-tower represented by
$(u',v',w')$. This sum also equals $g\circ T^h(y) -g(y)$.

\hfill $\square $

We can therefore give an alternate description of $K^0(X_B,T_B)$
which is more elegant than the description in Theorem 3.2. For the same reasons
as in the case of $\displaystyle\frac{\Z^{|V^{\mathcal
Q}_n|}}{B\Z^{|V^{\mathcal Q}_n|}}$, we see that the natural order
in $ \Z^{|V^{\mathcal Q}_n|}$  induces an order in
$\displaystyle\frac{\Z^{|V^{\mathcal
Q}_n|}}{\Delta\Z^{|V^{\mathcal Q}_n|}}$.

As we already observed, $A^{\mathcal Q}_n[\delta (\Z^{|V_n \times
V_n|})] \subseteq \delta (\Z^{|V_{n+1} \times V_{n+1}|})$.

\smallskip
\noindent{\bf 3.9. Theorem.} {\it For $B=(V,E,\geq)$ let
$(X_B,T_B)$ be defined as in (1.10). Write $(X,T)=(X_B,T_B)$.
Define the tripling $B^{\mathcal Q}=(V^{\mathcal Q},E^{\mathcal
Q},\geq)$ as in 2.2. Then the map induced between the inductive
limits of the two (horizontal) systems of ordered groups in the
following diagram is an isomorphism. }

$$\begin{matrix}
 \displaystyle\frac{\Z^{|V^{\mathcal Q}_0|}}{B\Z^{|V^{\mathcal Q}_0|}}
&\buildrel {A^{\mathcal Q}_0} \over \longrightarrow
&\displaystyle\frac{\Z^{|V^{\mathcal Q}_1|}}{B\Z^{|V^{\mathcal
Q}_1|}}& \buildrel {A^{\mathcal Q}_1} \over \longrightarrow &
\displaystyle\frac{\Z^{|V^{\mathcal Q}_2|}}{B\Z^{|V^{\mathcal
Q}_2|}} & \buildrel {A^{\mathcal Q}_2} \over \longrightarrow &
\displaystyle\frac{\Z^{|V^{\mathcal Q}_3|}}{B\Z^{|V^{\mathcal
Q}_3|}}&
\buildrel {A^{\mathcal Q}_3} \over \longrightarrow & \cdots      \\
&&&&&&&& \\
\rho_0\uparrow & & \rho_1\uparrow && \rho_2\uparrow && \rho_3\uparrow &&\\
&&&&&&&& \\
\displaystyle\frac{\Z^{|V^{\mathcal Q}_0|}}{\Delta\Z^{|V^{\mathcal
Q}_0|}}&\buildrel {A^{\mathcal Q}_0} \over \longrightarrow
&\displaystyle\frac{\Z^{|V^{\mathcal
Q}_1|}}{\Delta\Z^{|V^{\mathcal Q}_1|}}& \buildrel {A^{\mathcal
Q}_1} \over \longrightarrow & \displaystyle\frac{\Z^{|V^{\mathcal
Q}_2|}}{\Delta\Z^{|V^{\mathcal Q}_2|}} & \buildrel {A^{\mathcal
Q}_2} \over \longrightarrow & \displaystyle\frac{\Z^{|V^{\mathcal
Q}_3|}}{\Delta\Z^{|V^{\mathcal Q}_3|}}& \buildrel {A^{\mathcal
Q}_3} \over \longrightarrow & \cdots
\end{matrix}$$
{\it Furthermore, the two inductive limits are both isomorphic to
$K^0(X,T)$.}

[To avoid messy notation and display, we have hidden $G^+$ while
referring to the ordered group $(G,G^+$).]

\noindent{\bf Proof. }By lemma 3.6, $\Delta\Z^{|V^{\mathcal Q}_n|}
\subseteq B\Z^{|V^{\mathcal Q}_n|}$. By lemma 3.8, for
sufficiently large $K$, $A^{\mathcal Q}_{n+K-1}\circ \cdots
A^{\mathcal Q}_{n+1}\circ A^{\mathcal Q}_n ~(B\Z^{|V^{\mathcal
Q}_n|})\subseteq \Delta\Z^{|V^{\mathcal Q}_{n+K}|} $. Hence, the
induced map between the inductive limits is an isomorphism. That
$K^0(X,T)$ is isomorphic to the inductive limit of the top
horizontals  was already proved in theorem 3.2. Observe that $
\Z^{|V^{\mathcal Q}_0|}=\Z$ and $B\Z^{|V^{\mathcal
Q}_0|}=\Delta\Z^{|V^{\mathcal Q}_0|}=0 $. The image of $1\in
\Z^{|V^{\mathcal Q}_0|}$ in the inductive limit maps to the order
unit $u$ in $K^0(X,T)$ corresponding to the image of the constant
function $1 \in C(X,\Z)$.

\hfill $\square $

\smallskip

\noindent{\bf Remark. } Since it is known that $K^0(X,T)$ is
isomorphic to the $K_0$-group of the associated $C^*$-crossed
product $C(X)\rtimes_T \Z$, we see that as a corollary to theorem
3.9 , we can effectively compute $K_0(C(X_B)\rtimes_{T_B} \Z)$.

\smallskip

Finally, we should point out how these descriptions simplify
further for properly ordered Bratteli diagrams and yield the
isomorphism $K^0(X,T) \simeq K_0(V,E)$, (see 3.1), proved by
Hermann, Putnam and Skau [HPS, Theorem 5.4 and Corollary 6.3].

\smallskip

\noindent{\bf 3.10. }Let $(V,E,\geq)$ be a properly ordered Bratteli
diagram. Telescoping if necessary, assume that every level
$n+1$-tower traverses through at least two level $n$-towers. 
Telescoping further if necessary, (see [HPS, Proposition 2.8]), we can
assume that any two maximal edges of $E_n$ have the same source. Similarly,
we can assume that any two minimal edges  of $E_n$ have the same source.
For the rest of the paper we assume that these properties hold. Then
for any $\mathcal Q_{n+2}~$-tower $\mathcal S(u,v,w)$ the first
$\mathcal Q_n$-tower traversed by $\mathcal S(u,v,w)$ is
independent of $(u,v,w) \in V^{\mathcal Q}_{n+2}$. Thus one sees from the
definition of $B\Z^{|V^{\mathcal Q}_n|}$ that $A^{\mathcal
Q}_{n+1}\circ A^{\mathcal Q}_n (B\Z^{|V^{\mathcal Q}_n|}) = 0$. As
a consequence, the map induced between the inductive limits of the
top two horizontals in the following diagram is an isomorphism.

$$\begin{matrix}
 \Z
&\buildrel {A^{\mathcal Q}_0} \over \longrightarrow
&\displaystyle\frac{\Z^{|V^{\mathcal Q}_1|}}{B\Z^{|V^{\mathcal
Q}_1|}}& \buildrel {A^{\mathcal Q}_1} \over \longrightarrow &
\displaystyle\frac{\Z^{|V^{\mathcal Q}_2|}}{B\Z^{|V^{\mathcal
Q}_2|}} & \buildrel {A^{\mathcal Q}_2} \over \longrightarrow &
\cdots & \displaystyle\frac{\Z^{|V^{\mathcal
Q}_n|}}{B\Z^{|V^{\mathcal Q}_n|}}&
 \longrightarrow & \cdots      \\
&&&&&&&&&\\
1\uparrow & & \tau\uparrow && \tau\uparrow &&&\tau\uparrow &&\\
&&&&&&&&& \\
\Z & \buildrel {A^{\mathcal Q}_0} \over \longrightarrow
&\Z^{|V^{\mathcal Q}_1|}& \buildrel {A^{\mathcal Q}_1} \over
\longrightarrow & \Z^{|V^{\mathcal Q}_2|} & \buildrel {A^{\mathcal
Q}_2} \over \longrightarrow & \cdots & \Z^{|V^{\mathcal Q}_n|}&
\longrightarrow & \cdots \\
&&&&&&&&&\\
1\uparrow & & \pi^*\uparrow && \pi^*\uparrow &&&\pi^*\uparrow &&\\
&&&&&&&&& \\
\Z & \buildrel {A_0} \over \longrightarrow &\Z^{|V_1|}& \buildrel
{A_1} \over \longrightarrow & \Z^{|V_2|} & \buildrel {A_2} \over
\longrightarrow & \cdots & \Z^{|V_n|}& \longrightarrow & \cdots
\end{matrix}$$

\smallskip

The map $\pi^*:\Z^{|V_n|} \longrightarrow \Z^{|V^{\mathcal Q}_n|}$
is induced by the map $\pi:V^{\mathcal Q}_n \longrightarrow V_n$
given by $(u,v,w) \mapsto v$ and of course commutes with
multiplication by the respective incidence matrices (i.e,
$A^n,A^{\mathcal Q}_n$). 
Recall that after doing necessary telescoping we have arranged so
that the properly ordered Bratteli diagram $(V,E,\geq)$ has the properties
described in the beginning of 3.10.
Now let $e$ be an edge in $(V,E,\geq)$ with
range $v\in V_{n+1}$. Let $(u,v,w),(u',v,w') \in V^{\mathcal
Q}_{n+1}$. Let $\tilde{e},\tilde{e}'$ be the (unique) lifts of $e$
to $(V^{\mathcal Q},E^{\mathcal Q},\geq)$ with ranges
$(u,v,w),(u',v,w')$ respectively. Then, from the description in 2.2,
$\tilde{e},\tilde{e}'$ have the same sources in
$V^{\mathcal Q}_n$. From this it follows that for $\overline{p}
\in \Z^{|V^{\mathcal Q}_n|}, A^{\mathcal Q}_n
(\overline{p})(u,v,w)=A^{\mathcal Q}_n (\overline{p})(u',v,w')$;
in other words, $A^{\mathcal Q}_n (\overline{p}) \in \pi^*
(\Z^{|V_{n+1}| })$. Thus the map induced between the inductive
limits of the two bottom horizontals in the above diagram is also
an isomorphism.

\smallskip

\noindent{\bf 3.11. Specialization of Theorem 3.9 to
substitutional systems.} We recall the notation from 2.5. Let $(A,
\sigma)$ be a primitive aperiodic non-proper substitutional system. Let $B=(V,E,\geq )$ be the
stationary ordered Bratteli diagram associated to $(A, \sigma)$.
Define $ A^{\mathcal Q} = \{(a,b,c) \in A \times A \times A \mid
abc \text{ occurs as a subword of } \sigma^n (d) \text{ for some
}d \in A \text{ and some } n \} $. Define
$$\sigma^{\mathcal Q}: A^{\mathcal Q} \to (A^{\mathcal Q})^+$$
by $\sigma^{\mathcal Q} [(a,b,c)] =(a_m, b_1, b_2) \cdot
(b_1,b_2,b_3) \cdots (b_{n-2}, b_{n-1}, b_n) \cdot (b_{n-1}, b_n,
c_1)$,  where $\sigma (b) =b_1 \cdot b_2 \cdots b_n$, and $a_m$ is
the last alphabet in $\sigma (a)$, while $c_1$ is the first
alphabet in $\sigma (c)$. For a function $\varphi :A\times A
\longrightarrow \Z$ define $\delta (\varphi) \in \Z^{|A^{\mathcal
Q}|}$ by $\delta (\varphi)(a,b,c) = \varphi (b,c) - \varphi
(a,b)$. Let $\Delta \Z^{|A^{\mathcal Q}|}$ be the subgroup $\delta
(\Z^{|A \times A|})$  of $\Z^{|A^{\mathcal Q}|}$. Suppose
$\overline{p},\overline{q} \in \Z^{|A^{\mathcal Q}|}$ and take
values in $\Z^+$ and further that
$\overline{p}=-\overline{q}~${\bf mod } $\Delta \Z^{|A^{\mathcal
Q}|}$. Then $\overline{p}=\overline{q}=0$. The natural order in
$\Z^{|A^{\mathcal Q}|}$ induces an order in $\Z^{|A^{\mathcal
Q}|}/\Delta \Z^{|A^{\mathcal Q}|}$ making it an ordered group.

\smallskip

\paragraph{} Let $\beta^{\mathcal Q}: \Z^{|A^{\mathcal
Q}|}\longrightarrow \Z^{|A^{\mathcal Q}|}$ be given by matrix
multiplication by the incidence matrix of the substitution
$\sigma^{\mathcal Q}$. For a function $\varphi :A \times A
\longrightarrow \Z$ define $\varphi':A \times A \longrightarrow
\Z$ by $\varphi' (u',v')= \varphi (u,v)$ where $u$(resp.$v$) is
the last(resp.first) alphabet in the substitution $\sigma
(u')$(resp.$\sigma (v')$). Then $\beta^{\mathcal
Q}(\delta(\varphi)) =\delta(\varphi')$; thus, $\beta^{\mathcal
Q}[\delta (\Z^{|A\times A|})] \subseteq \delta (\Z^{|A \times
A|})$. Hence, $\beta^{\mathcal Q}$ induces a homomorphism of
ordered groups
$$\displaystyle \frac{\Z^{|A^{\mathcal Q}|}}{\Delta \Z^{|A^{\mathcal
Q}|}} \longrightarrow \displaystyle \frac{\Z^{|A^{\mathcal
Q}|}}{\Delta \Z^{|A^{\mathcal Q}|}}$$ still denoted by
$\beta^{\mathcal Q}$.

\smallskip

From theorem 3.9 one deduces immediately

\noindent{\bf Theorem 3.12. }{\it The inductive limit of the
system of ordered groups
$$\begin{matrix}

\displaystyle\frac{\Z^{|A^{\mathcal Q}|}}{\Delta\Z^{|A^{\mathcal
Q}|}}&\buildrel {\beta^{\mathcal Q}} \over \longrightarrow
&\displaystyle\frac{\Z^{|A^{\mathcal Q}|}}{\Delta\Z^{|A^{\mathcal
Q}|}}& \buildrel {\beta^{\mathcal Q}} \over \longrightarrow &
\displaystyle\frac{\Z^{|A^{\mathcal Q}|}}{\Delta\Z^{|A^{\mathcal
Q}|}} & \buildrel {\beta^{\mathcal Q}} \over \longrightarrow &
\displaystyle\frac{\Z^{|A^{\mathcal Q}|}}{\Delta\Z^{|A^{\mathcal
Q}|}}& \buildrel {\beta^{\mathcal Q}} \over \longrightarrow &
\cdots
\end{matrix}$$
is isomorphic to the dimension group $K^0(X_{\sigma},T_{\sigma})$ of the substitution system
associated to $(A,\sigma)$.}

\bigskip

\bigskip

\centerline{\bf References}

\bigskip

\begin{enumerate}
\item[[DHS]] {\sc F. Durand, B. Host and C. Skau}. Substitutional
dynamical systems, Bratteli diagrams and dimension groups, {\it
Ergodic Th. and Dynam. Sys.} {\bf 19}, (1999) 953-993.

\item[[EP]] {\sc A. El Kacimi, R. Parthasarathy}. Skew-product for
group-valued edge labellings  of Bratteli diagrams, {\it
arXiv.math.DS/0506304}.

\item[[HPS]] {\sc R.H. Herman, I.F. Putnam and C.F. Skau}. Ordered
Bratteli diagrams, dimension groups and topological dynamics, {\it
Internat. J. Math.} {\bf 3}, (1992) 827-864.
\end{enumerate}
\end{document}